\input epsf

\baselineskip = 0.2in
\centerline {\bf  KNOTTED CONTRACIBLE 4-MANIFOLDS IN THE 4-SPHERE}

\vskip 0.3truein
\centerline {\bf by}
\vskip 0.3truein  

\centerline {\bf  W.B.Raymond Lickorish}

\vskip 0.5truein

Abstract \quad  Examples are given to show that some compact contractible 4-manifolds
can be knotted in the 4-sphere.
It is then proved that any finitely  presented perfect group with a balanced presentation is a knot
group for an embedding of some contractible 4-manifold in $S^4$.
\vskip 0.5truein

\noindent {\bf 1. Introduction}

\vskip 0.3truein

A construction will here be described that can produce a compact contractible 4-manifold 
$M$ embedded piecewise linearly (or smoothly) in
$S^4$ with  the fundamental group of its complement being non-trivial.  Then, another embedding of
$M$ in $S^4$ will be produced which has  simply connected complement.  Several examples of this will be
given.  Of course, the construction
  emphasises that contractible spaces do not behave entirely as do single
points.  
It is important to note
that these embeddings are piecewise linear or smooth; they are certainly not wild.  The famous
construction of the Alexander wild horned sphere gives a wild embedding of a 3-ball in $S^3$ that

has its complement not simply connected.  However the boundary of a contractible compact 
3-manifold is just a 2-sphere,  so by the  piecewise linear 3-dimensional
Sch\"onflies theorem, if such a manifold can be embedded {\it piecewise linearly} in
$S^3$, each of the manifold and its complement must be a 3-ball.

 Recall the general definition of knotting, when all maps and spaces are in the
piecewise linear category:  A polyhedron
$X$ knots in a polyhedron
$Y$ if there are two embeddings, $e_0$ and $e_1$ of $X$ in $Y$, that are homotopic but not ambient isotopic.  The
embeddings are ambient isotopic if there exist homeomorphisms $F_t : Y \rightarrow Y$, for each $t \in
[0,1]$, such that $(y, t) \mapsto (F_t(y), t)$  defines a piecewise linear homeomorphism
from $ Y \times [0,1] $ to itself, $F_0$ is the identity and $F_1e_0 = e_1$. Thus to
be ambient isotopic the complements of the images of the two embeddings must certainly be
homeomorphic.  The knotting phenomenon explores the possibility of moving between embeddings
along a path of embeddings as opposed to moving along a path of maps.  The examples given here are of
contractible 4-manifolds that can knot in
$S^4$ for, just as in classical knot theory, the fundamental group of  complements is used to show
embeddings are not ambient isotopic.  Examples of knots usually rely on the entwining of some non-trivial
cycle, but here there is none. In fact, in higher dimensions, if
$X$ and
$Y$ are piecewise linear manifolds  with
${\hbox {dim}\,}Y -  {\hbox {dim}\,}X \geq 3$,  there are theorems of Hudson [4] that assert that
there is no knotting of $X$ in $Y$ provided these spaces are sufficiently highly connected.

It should be noted that  when $M$ is a contractible 4-manifold, piecewise linearly contained in
$S^4$, the Alexander duality theorem implies that  $S^4 - M$ has the same homology as a point.  Thus
$\pi_1(S^4 - M)$ is a perfect group in contrast to the situation of classical knot theory.  It will be
proved in Theorem 3 that for any  perfect group with a balanced presentation (that is, a presentation with the same finite number of generators as relations), there are embeddings $e_0$ and $e_1$ of some a contractible 4-manifold  $M$ into
$S^4$ so that 
$\pi_1(S^4 - e_1M)$ is the given group and $\pi_1(S^4 - e_0M)$ is trivial.

The author is grateful  to Simon Norton for making a helpful remark and to  Charles Livingston for a correction.

\bigskip

\vskip 0.2truein

\noindent{\bf 2. Examples of the embedding construction}
\vskip 0.3truein
The  theorem that now follows is really an example  describing the main simple idea of the construction of
this paper. The second theorem amplifies it to more general circumstances.

\proclaim {Theorem 1}.  There are two piecewise linear (or smooth) embeddings, $e_0$ and $e_1$, of a
certain compact contractible 4-manifold
$M$ into $S^4$ such that $\pi_1(S^4 - e_1M)$ is non-trivial and $S^4 - e_0M$ is contractible.

\noindent  {\it Proof} \quad Firstly, construct a compact 4-manifold $X$ by adding three 1-handles and three
2-handles onto a 4-ball in the following way.  The handles are to be chosen so that $\pi_1(X)$ has
the presentation 
$$\langle a, b, c : \ b^{-1}c^{-2}bc^3, \  c^{-1}a^{-2}ca^3,  \  a^{-1}b^{-2}ab^3 \rangle    \ ,$$  
where based loops encircling the three 1-handles represent $a$, $b$ and $c$, and the attaching
circles of the three 2-handles give the three relators.  This situation is shown   in
Figure 1 in the notation common in considerations of the `Kirby calculus' (see  [2] for example).

\epsfxsize=2in
\centerline
{\epsfbox  {mainidea}}
\centerline{Figure 1}
  
\vskip 0.1in

The   diagram of Figure 1 shows curves in the 3-sphere, the boundary of the 4-ball.  Open regular
neighbourhoods, of three standard disjoint discs in the 4-ball, are to be removed from the 4-ball to
create a ball with the three 1-handles added.
The boundaries of these discs are the  circles,
decorated with dots, labeled $a$,
$b$ and $c$.   That
this is in order can be checked as follows. A ball with   1-handles added can be changed back to a
ball by adding  2-handles to cancel the 1-handles; removing those 2-handles consists of removing
neighbourhoods of the discs that are the  co-cores of the 2-handles.  Thus a 4-ball, with
1-handles added, is the same as a 4-ball from which standard 2-handles have been removed.  A
1-handle can be regarded as
$D^1
\times D^3$ with $\partial D^1 \times \star
$ being the attaching sphere and $\star \times \partial D^3$ being the belt sphere (where each
$\star$ is a base point).  In Figure 1 a belt sphere consists of the union of a disc spanning a
dotted circle, less a regular neighbourhood of that circle, and a  disc in the boundary of the
2-handle that has been removed.  Meridians encircling the three dotted circles represent
generators, to be called
$a$,
$b$ and
$c$, of the fundamental group of the ball with 1-handles, and a based closed curve represents a
word, in $a$,
$b$ and $c$, corresponding to its signed intersections with the three belt spheres.  In this way the
curves shown, labeled $\alpha$,
$\beta $ and $\gamma$, represent  $b^{-1}c^{-2}bc^3$, $c^{-1}a^{-2}ca^3$ and $a^{-1}b^{-2}ab^3$.
Thus adding 2-handles, with these curves as attaching spheres (choose the zero framings),
gives the 4-manifold
$X$ with the required presentation for $\pi_1(X)$.  It has been shown by Rapaport [6] that this is
the presentation of a {\it non-trivial group}.  There are, of course, very many ways that attaching
curves can be chosen for  the 2-handles in order to achieve this presentation (and the choice will be
explored further in Theorem 3), but the one shown is about the simplest and is the one that will now be
considered.  The situation is shown schematically in Figure 2.

\epsfxsize=1.70in
\centerline
{\epsfbox{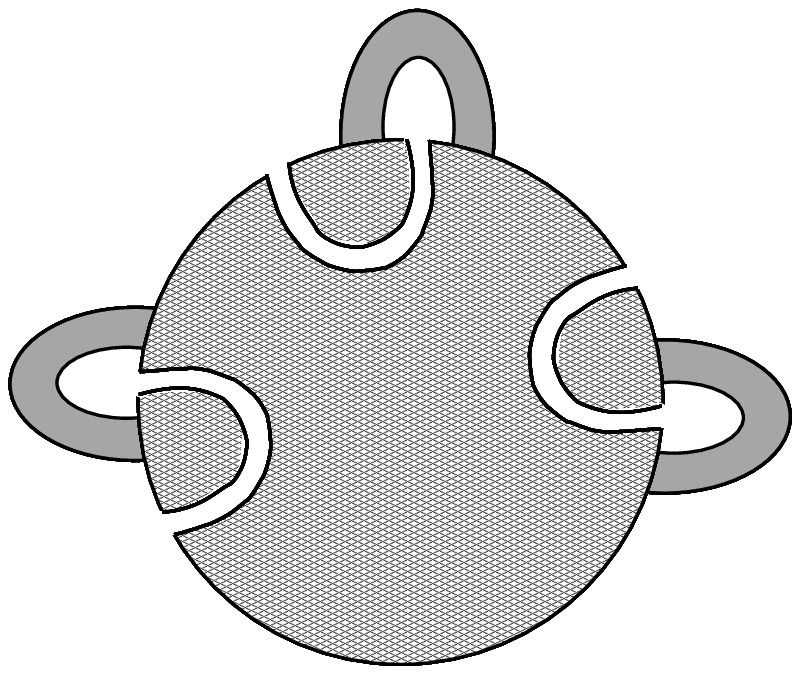}}
\centerline{Figure 2}
  
 \vskip 0.1in

The manifold $X$ is a 4-ball, $B_1$ say,  with 2-handles removed and 2-handles added.  Regard this
4-ball as being contained in $S^4$ and consider the complementary 4-ball $B_2$.  The 2-handles
removed from   $B_1$  can be thought of as added to $B_2$.  The other 2-handles, that were added to 
$B_1$, were added with zero framing along  unknotted, unlinked curves (labeled $\alpha$,
$\beta $ and $\gamma$), so they can be regarded as standard 2-handles removed from $B_2$.  Thus the
closure of $S^4 - X$ is the 4-ball $B_2$ with three 2-handles removed (creating added 1-handles) and
three 2-handles added and {\it this} is to be the required 4-manifold
$M$.  The situation for $M$ is again represented by Figure 1, except that now the dots should be
removed from the curves labeled 
$a$, $b$ and
$c$ and placed on those labeled  
$\alpha$, $\beta $ and $\gamma$.  However, the  $\alpha$-curve bounds a disc that meets only the
$c$-curve; there are similar discs for the $\beta$- and $\gamma$-curves.  
The words in $\alpha$, $\beta $ and $\gamma$, coming from the intersections of the $a$, $b$ and
$c$ curves with these discs, make it clear that 
$\pi_1(M)$ is presented by 
$$\langle \alpha, \beta, \gamma \ : \ \alpha^3\alpha^{-2}, \ \beta^3\beta^{-2}, \
\gamma^3\gamma^{-2} \rangle $$ which, very obviously, presents the trivial group.  Thus $M$ is
simply connected.  A count of the handles shows that the Euler characteristic of $M$ is 1, hence
$H_2(M) = 0$.  Furthermore  $H_r(M) = 0$ for $r > 2$, as   there are no $r$-handles for $r > 2$, and
so, by the Hurewicz isomorphism theorem, $M$ has all homotopy groups trivial and hence is
contractible.  Note that
$\pi_1(\partial M)
\neq
\{ 1
\}$, as otherwise $\pi_1(X)
\cong \pi_1(S^4)$ by the Van Kampen theorem.  Hence $M$ is not a 4-ball.  The inclusion of $M$, as
so defined in $S^4$, is the embedding
$e_1$.  

The above {\it presentation} given for  $\pi_1(M)$ coming from the handle structure of $M$ is almost
trivial.  It certainly reduces to the trivial presentation by Andrews-Curtis moves (see below). 
Any such
$M$ has the property that
$M \times [0,1] \cong B^5$ where   $B^5$ is a 5-ball.  To show that in this instance, it is  necessary
only to realise that $M \times [0,1]$  has the same handle structure as does $M$. The extra
dimension means that, when a 2-handle is attached (to the boundary of a 5-manifold) only the
homotopy class of the attaching map is significant (a homotopy of attaching circles can be changed
to an isotopy by using the fourth dimension to  prevent the circles from crossing each other).  Now
let
$e_0$ be the inclusion of
$M
\times
\{ 0 \}$ in the 4-sphere
$\partial (M
\times [0,1])$.  The complement of $M \times \{ 0 \}$  in this sphere is 
$(\partial M \times [0,1]) \cup M \times \{ 1 \}$ and this is just another copy of the contractible
manifold
$M$. 
 
\vskip 0.3truein

For a second example consider 
$\langle a, b : \ ab^2ab^{-1}, \ a^4ba^{-1}b  \rangle  $, a presentation of the perfect group $G$ of
120 elements that is the fundamental group of the Poincar\'e homology 3-sphere.  The method
of the proof of Theorem 1 constructs a 4-manifold $X \subset S^4$ with $\pi_1(X) \cong G$ and with
the fundamental group of the corresponding $M$ presented by 
$\langle \alpha , \beta : \alpha ^2\beta ^3, \ \alpha ^{-1}\beta ^{-2}   \rangle  $.  Again this
easily reduces to the trivial presentation by Andrews-Curtis moves so that $M \times I$ is a 5-ball.

The  above construction works easily for the presentations 
$$\langle a_1,a_2, \dots, a_m : r_1,r_2, \dots, r_m 
\rangle $$ for every $m \geq 4$ when $r_i = a_i^{-1}a_{i+1}^{-1}a_ia_{i+1}^2$ for $i = 1, 2, \dots
m$ modulo $m$, and also when $r_i = a_i^{-1}a_{i+1}^{-2}a_ia_{i+1}^3$.  These  are known to
be presentations of infinite groups (see [3] and [5]).

Note that $\langle a, b : a^{-1}b^{-2}ab^3, b^{-1}a^{-2}ba^3 \rangle $ is a presentation of the
trivial group.  If $M$ is constructed from this presentation for $X$ it is not clear whether the
embedding of $M$ is in any sense knotted.  
\vskip 0.3truein

\noindent {\bf 3. A reminder of  Andrews-Curtis moves}
\vskip 0.3truein
The above proof makes a brief mention of the Andrews-Curtis moves. 
These moves are elementary changes that can be made to a group presentation 
that do not alter the group that is presented.  The moves are also called `extended Nielsen
transformations' in [1], they are called `$Q$-transformations' in [6] and they are sometimes
also called `Markov operations'.  The permitted changes to a presentation
$\langle a_1,a_2, \dots, a_m : r_1,r_2, \dots, r_n  \rangle $  are the following moves 
and  the inverses of these moves.

\begingroup
\parindent = 0.5truein
\item {(i)}  Change $r_i$ to $r_ia_ja_j^{-1}$ or $r_ia_j^{-1}a_j$.
\item {(ii)} Change $r_i$ to a cyclic permutation of $r_i$.
\item {(iii)} Change $r_i$ to $r_i^{-1}$.
\item {(iv)} Change $r_i$ to $r_ir_j$ where $j \neq i$.
\item {(v)}  Add a new generator $a_{m+1}$ and a new relator $a_{m+1}w$ where $w$ is a word in
$a_1,a_2,\dots, a_m$.
\endgroup

These are precisely the moves that can easily be imitated on a 5-manifold comprised of 0-handles,
1-handles and 2-handles only. If the handles of such a 5-manifold correspond to a presentation of
the trivial group that can be reduced to the trivial presentation (that is, the empty
presentation) by the above moves and their inverses, then it is shown in [1] that the manifold is
the 5-ball.  It is this result that is used in the above proof.  The Andrews-Curtis conjecture [1]
is that {\it any} presentation of the trivial group be reducible to the trivial presentation by the
above moves and their inverses.  This is popularly thought to be false,
$\langle a, b : a^{-1}b^{-2}ab^3, b^{-1}a^{-2}ba^3
\rangle $ being one of many proposed counter-examples.  The truth of the Andrews-Curtis conjecture 
would imply the truth
of  another conjecture that asserts that any 5-dimensional regular neighbourhood of a contractible
2-complex be a 5-ball (such a neighbourhood is known to be unique).

\vskip 0.3truein

\noindent {\bf 4. Arbitrary finitely presented perfect groups}

\vskip 0.3truein

A few simple remarks lead up to an elementary, but possibly surprising, little lemma about finitely
presented perfect groups and presentations of the trivial group. Suppose that an {\it abelian} group
$E$ is freely generated as an {\it abelian} group (with additive notation) by the generators $e_1, e_2,
\dots, e_m$. The quotient of $E$ by the subgroup generated by the $n$ elements 
$\{ \ \sum ^m_{j=1}a_{i j}e_j : i = 1,2,\dots,n \ \} $ is said to be presented by the $n \times m$
integer matrix 
$A = \{ a_{i j} \}$.
 When $n = m$ the quotient is the trivial group if and only if $A$ is unimodular, that is, $\det A = \pm 1$. 
A presentation $P$ of any group (in multiplicative notation) leads at once to a
presentation of the abelianisation of that group, by just deciding that all symbols
commute.  It is then sensible in each relator to assemble together all occurrences of a
generator and its inverse, cancelling where possible,  to obtain from the resulting
exponents in each relator  a presentation matrix $A$ of the abelianisation of the
group.   The following Lemma considers such things in the reverse order, showing that, if
$A$ presents the trivial abelian group, then $P$ can be chosen to present the trivial
group.

\proclaim { Lemma 2}. Suppose that $A$ is a unimodular $n \times n$ matrix of integers. 
Then there exists a presentation $P$ of the trivial group that has $A$ as its
abelianised presentation matrix.  Furthermore, $P$ is equivalent to the trivial
presentation by Andrews-Curtis moves.

\noindent {\it Proof} \quad  Starting from the identity $n \times n$ matrix, the unimodular matrix $A$ can
be created by a sequence of row operations in which either a row is multiplied by $-1$ or
a row is added to another row.  These moves can be mimicked by changes to a
presentation 
$\langle a_1, a_2, \dots, a_n  :  r_1, r_2, \dots, r_n \rangle$
 of the trivial group,
 where initially $r_i = a_i$ for each
$i$.  If the $i$th row of the matrix is multiplied by $-1$, change $r_i$ to $r_i^{-1}$; if row $i$ is
added to row $j$ then change $r_j$ to $r_jr_i$.  At each stage the presentation is of the trivial group
and at each stage the matrix is the {\it corresponding} presentation matrix of the  abelianised group. 
Of course the moves used on the presentation are all Andrews-Curtis moves.
\vskip 0.3truein

\proclaim {Theorem 3}.  Let $G$ be any  perfect group having a finite balanced presentation.  Then there is a compact
contractible 4-manifold $M$ contained in
$S^4$ such that 
$\pi_1(S^4 - M) \cong G$ and $M \times I$ is a 5-ball (so that, if $G$ is non-trivial, there is a
distinct second embedding of
$M$ in
$S^4$ having contractible complement).

\noindent  {\it Proof}\quad  Let $\langle a_1, a_2, \dots, a_n  :  r_1, r_2, \dots, r_n \rangle$ be a presentation $P$ of $G$.  The
construction proceeds as in the proof of Theorem 1.  Remove from the 4-ball $B^4$ neighbourhoods of $n$ standard 
disjoint spanning discs  to create a ball with $n$ 1-handles added.  
The boundaries of the discs form a set 
of $n$ unlinked simple closed (`dotted') curves in $\partial B^4 = S^3$, which are labelled $a_1, a_2, \dots, a_n$. 
In the following way  construct, as the boundaries of disjoint discs $D_1, D_2, \dots, D_n$  contained in
$S^3$,  simple closed curves $\alpha_1, \alpha_2, \dots, \alpha_n$, 
corresponding to 
$r_1, r_2, \dots, r_n$, which are to be the attaching circles for $n$   2-handles.   
Begin with small,
disjoint, oriented discs $\Delta_1, \Delta_2, \dots, \Delta_n$ in the complement of $a_1 \cup a_2 \cup \dots \cup
a_n$.  For each letter $a_i^{\pm 1}$  
  in the word $r_1$ take a small meridian disc of the curve $a_i$,
oriented according to the exponent on the letter, and  to construct $D_1$, join the boundaries of
these meridian discs by thin bands to the boundary of $\Delta_1$, in the order around $\partial \Delta_1$ specified
by
$r_1$.  The discs $D_2, \dots, D_n$ 
are constructed similarly from 
$r_2, \dots, r_n$ and there is no difficulty in
ensuring that the $D_i$ are embedded and mutually disjoint.

As in Theorem 1, form a 4-manifold $X \subset S^4$ by adding $n$ 2-handles with zero framing along
$\alpha_1, \alpha_2, \dots, \alpha_n$
to the ball with $n$ 1-handles.  Then $\pi_1(X) \cong G$.  The key point to note now is that the
meridian discs described above (for all the $D_i$  together) can be taken {\it in any order}
around $a_i$.  Different choices of order probably give different manifolds $X$ and $M$, where
again $M$ is the closure of  $S^3 - X$. Also note that the given presentation of $G$ can be amended by the insertion
of any number of copies of $a_ia_i^{-1}$, for any $i$, into any $r_j$ without changing $G$ nor the presentation
matrix $A$ of the (trivial) abelianisation of $G$ coming from that presentation. 
Now $\pi_1(M)$ has a
presentation $\Pi$ of the form $\langle \alpha_1, \alpha_2, \dots, \alpha_n : \rho_1, \rho_2, \dots, \rho_n \rangle$
where the relators record in order the occurrence of the meridian  discs  around
$a_1, a_2, \dots, a_n$ (each signed intersection of $a_i$ with a meridian disc contained in
$D_j$ producing an $\alpha_j^{\pm 1}$ entry  in $\rho_i$).  The abelian presentation matrix coming from $\Pi$ is the
transpose of $A$; it is certainly unimodular. Thus using Lemma 2, the ordering along the $a_i$ of those meridional
discs making up each $D_j$ can be chosen, after inserting any necessary pairs of discs  corresponding to
$a_ia_i^{-1}$, so that, with respect to the new choice, $\Pi$ becomes a presentation of the
trivial group.  Again from Lemma 2, $\Pi$ is equivalent by Andrews-Curtis moves to the trivial presentation and so
$M \times I$ is a 5-ball.

\vskip 0.3truein

\noindent {\bf References}
\vskip 0.1truein  

\noindent [1] \quad J.J.Andrews and M.L.Curtis. \   Free groups and handlebodies, Proc. Amer. Math. Soc. 16 (1965)
\hbox{192--195.}
\vskip 0.1truein
\noindent [2] \quad E.C\'esar de S\'a. A link calculus for $4$-manifolds,  Lecture Notes in Math. 722,  Springer, Berlin (1979)\hbox {16 -- 30.} 
\vskip 0.1truein
\noindent [3] \quad G.Higman.  A finitely generated infinite simple group, J. London Math. Soc. 26 (1951) 61--64.
\vskip 0.1truein
\noindent [4]  \quad J.F.P.Hudson.  {\it Piecewise linear topology},Benjamin, New York (1969).
\vskip 0.1truein
\noindent [5]\quad B.H.Neuman.   An essay on free products of groups with amalgamation, Philos. Trans. Roy. Soc. London A246 (1954) 503--554.
\vskip 0.1truein
\noindent [6]\quad E.S.Rapaport.  Groups of order 1, Proc. Amer. Math. Soc. 15 (1964) 828--833.

\vskip 0.3truein

\noindent University of Cambridge,
\hfil\break 
\noindent Department of Pure Mathematics and Mathematical Statistics
\hfil\break
\noindent Centre for Mathematical Sciences 
\hfil\break
\noindent Wilberforce Road, Cambridge CB3 0WB,  U.K.

\noindent {\bf email}  wbrl@dpmms.cam.ac.uk

\end